\documentclass[12pt]{amsart}
\usepackage[russian,english]{babel}
\usepackage{amsmath,amssymb,amsthm,amscd}
\usepackage[cp866]{inputenc}
\usepackage{graphicx}
\usepackage{float}
\usepackage{hyperref}

\theoremstyle{plain}
\newtheorem{theorem}{Theorem}
\newtheorem{lemma}[theorem]{Lemma}

\theoremstyle{definition}
\newtheorem{definition}{Definition}

\theoremstyle{plain}
\newtoks\thehProclaim
\newtheorem*{Proclaim}{\the\thehProclaim}

\theoremstyle{definition}
\newtoks{\thehRemark}
\newtheorem*{Remark}{\the\thehRemark}

\renewcommand{\leq}{\leqslant}
\renewcommand{\geq}{\geqslant}

\begin{document}

	\title{Minimal quasisections of circle bundles}

\author{Ivan Kiselev, Gaiane Panina, Timur Shamazov}
	
\address{G. Panina: St. Petersburg Department of Steklov Mathematical Institute; gaiane-panina@rambler.ru; I. Kiselev: St. Petersburg State University,kiiv2001@mail.ru; Timur Shamazov:  St. Petersburg Department of Steklov Mathematical Institute; shamazovt@inbox.ru}

	\subjclass[2000]{}
	

	\begin{abstract}
Let $E\xrightarrow{\pi}B$ be an oriented circle bundle over an oriented
closed surface, and let $\mathcal E$ denote its Euler number. A
\textit{quasisection} is a smooth surface $Q$, possibly with boundary,
together with a generic smooth map $q:Q\to E$ such that
$\pi\circ q(Q)=B$.

The \textit{local formula} expresses the Euler number
$\mathcal E$ of a bundle as a sum of the weights of the \textit{singular vertices} of a
quasisection.

In the present paper determine the minimal number $\mathcal N$ of singular vertices for
closed and for bordered quasisections for a given circle bundle. For
bordered quasisections,
\[
 \mathcal N_{\mathrm b}(E\to B)
 =\left\lceil\frac{3|\mathcal E|}{2}\right\rceil.
\]

\[
 \mathcal N_{\mathrm{cl}}(E\to B)=
 \begin{cases}

  \dfrac{3|\mathcal E|}{2},&\mathcal E\text{ is even},\\
  4|\mathcal E|,&\mathcal E\text{ is odd}.
 \end{cases}
\]

We also provide examples of quasisections that realize these minimal values.

	\end{abstract}

	\maketitle

	\section{Introduction}\label{SecIntro}

	Let $E\xrightarrow{\pi}B$ be an oriented circle bundle, that is, an
	oriented locally trivial fiber bundle with fiber $S^1$ over a closed oriented surface \cite{Chern}.
	Such bundles are
	classified by their Euler classes, or, equivalently in this setting,
	by their Euler numbers.
	
\medskip

	\begin{definition} \label{DefNumberSing}
		A \emph{closed quasisection} of a circle bundle 
		$E \xrightarrow[\text{}]{\pi} B$ is a pair
		$\mathcal{Q}=(Q,q)$, where $Q$ is a closed smooth surface and
		$q:Q\rightarrow E$ is a smooth map such that
		$\pi\circ q(Q)=B$. A \emph{bordered quasisection} is defined in
		the same way, for a compact smooth surface $Q$ possibly with
		boundary. In both cases, $Q$ is not required to be connected.
		
	\end{definition}
	Throughout the paper, the maps $q$ and $\pi\circ q$ are assumed to be
	generic, or stable, in the sense of Whitney singularity theory
	\cite{AVGZ}. Thus the image of $q$ may have curves of transverse self-intersection,
	isolated triple points, and Whitney umbrellas. At regular points of
	$q$, the only singularities of $\pi\circ q$ are folds and pleats. The
	additional genericity assumptions used below are listed in
	Section~\ref{SecSing}.

Given a quasisection, Theorems~\ref{ThmMain} and~\ref{ThmBorderedLCF}, taken from \cite{PanShamTur}, express the Euler number of the bundle as the
sum of the weights of certain singularities of the quasisection. This singularities are
called \textit{singular vertices}, see section~\ref{SecSing}--\ref{SecWeightsClosed} for definition.

\begin{figure}[!b]
  \centering
  \includegraphics[width=0.72\textwidth]{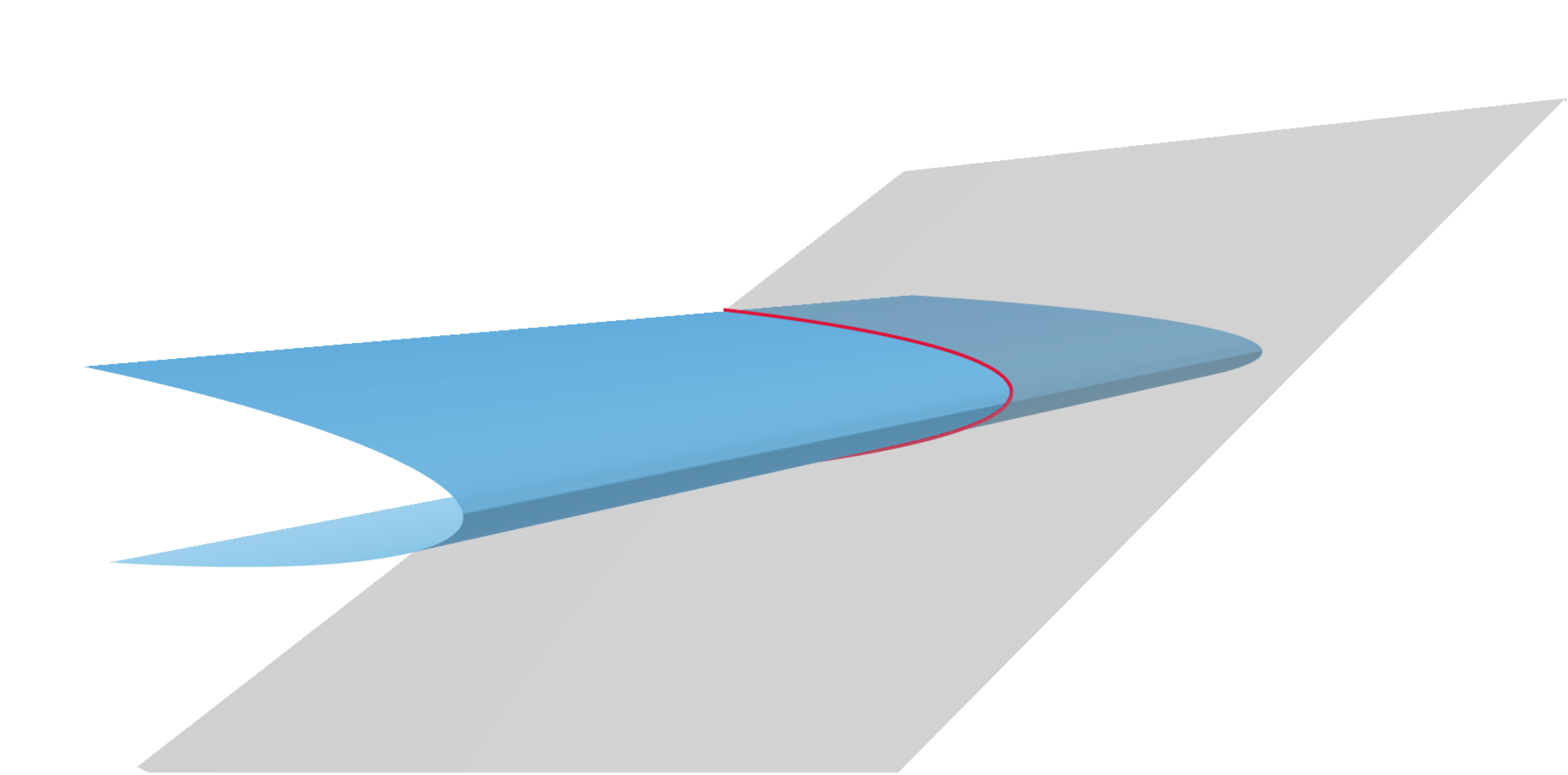}
  \caption{A regular sheet crosses a fold.}\label{FigCross}
\end{figure}

\begin{definition}
	The \textit{number of singularities} $N(\mathcal{Q})$ is the number of
	all singular vertices of $\mathcal{Q}$.
	We define two numbers:
	\begin{align*}
	\mathcal{N}_{\mathrm{cl}}(E\rightarrow B)
	&=\min\{N(\mathcal{Q}):\mathcal{Q}\text{ is a closed quasisection of } E\rightarrow B\},\\
	\mathcal{N}_{\mathrm{b}}(E\rightarrow B)
	&=\min\{N(\mathcal{Q}):\mathcal{Q}\text{ is a bordered quasisection with }\partial Q\ne\varnothing\}.
	\end{align*}
	Clearly, $\mathcal{N}_{\mathrm{cl}}$ and $\mathcal{N}_{\mathrm{b}}$ are functions of $\mathcal{E}$. 
	
	A closed (respectively, bordered) quasisection is \textit{minimal in
	its class} if its number of singular vertices equals
	$\mathcal{N}_{\mathrm{cl}}(E\rightarrow B)$ (respectively,
	$\mathcal{N}_{\mathrm{b}}(E\rightarrow B)$).
	\end{definition}

In this paper, we establish the values of $\mathcal{N}_{\mathrm{b}}$ and $\mathcal{N}_{\mathrm{cl}}$ and provide examples of quasisections that realize these minimal values.

\begin{theorem}\label{ThmMinimal}
For a circle bundle $E\rightarrow B$ with Euler number $\mathcal{E}$,
\begin{enumerate}
\renewcommand{\labelenumi}{(\arabic{enumi})}
\item\label{ThmMainBordered}
\begin{minipage}[c]{\linewidth}
\[
 \mathcal{N}_{\mathrm{b}}(E\rightarrow B)
 =\left\lceil\frac{3|\mathcal{E}|}{2}\right\rceil,
\]
\end{minipage}

\item\label{ThmMainClosed}
\begin{minipage}[c]{\linewidth}
\[
 \mathcal{N}_{\mathrm{cl}}(E\rightarrow B)=
 \begin{cases}

  \dfrac{3|\mathcal{E}|}{2}, & \mathcal{E}\text{ is even},\\
  4|\mathcal{E}|, & \mathcal{E}\text{ is odd}.
 \end{cases}
\]
\end{minipage}
\end{enumerate}
\end{theorem}

\section{Preliminaries}\label{SecPrel}

	\medskip

	We recall the main constructions and definitions, referring the reader to
\cite{PanShamTur} for further details, figures,
and examples.
	
	\subsection{Singular vertices of a closed quasisection}\label{SecSing}

We first recall the constructions for closed quasisections: the
singular curve, the portraits of singular vertices, their weights, and
the local formula. Throughout Sections~\ref{SecSing}--\ref{SecWeightsClosed},
the quasisection is assumed to be closed. In
Section~\ref{SecBordered} we then describe precisely what changes when
the quasisection is allowed to have boundary.

Let $\mathcal{Q}=(Q,q)$ be a closed quasisection of a circle bundle
$E\rightarrow B$.

The \textit{singular curve} of $\mathcal Q$ is the closure in $B$ of
the union of the following sets:
\begin{enumerate}
 \item the image under $\pi$ of the self-intersection locus of $q(Q)$;
 \item the image under $\pi\circ q$ of the fold locus.
\end{enumerate}
We refer to the second set as the \textit{projection of the folds}.
In particular, the images under $\pi\circ q$ of Whitney umbrellas and pleat points belong
to the singular curve; we call them the \textit{projections of Whitney
umbrellas} and the \textit{projections of pleats}, respectively.

\medskip

	A \textit{singular vertex} is a point of one of the following types:
	\begin{enumerate}
		\item A self-intersection of the singular curve.
		\item The projection of a Whitney umbrella.
		\item The projection of a pleat.
	\end{enumerate}

We assume that all singular vertices are \textit{stable}, that is, that they
persist under sufficiently small perturbations of $q$. The precise genericity conditions are listed in \cite{PanShamTur}.

	\subsection{Portraits of singular vertices of a closed quasisection}\label{SecPortraits}
	
	Let $\mathcal Q=(Q,q)$ be a closed quasisection and let $x\in B$.
Choose a small positively oriented circle $C_x$ around $x$. Equip a
neighborhood of $C_x$ and its preimage with product metrics for which
$\pi$ is locally an orthogonal projection and all fibers have the same
length. Then $\pi^{-1}(C_x)$ is a flat torus. Cut this torus along a
circle transverse to the fibers and depict it as an annulus whose inner
and outer boundary circles are identified. Orient the fibers inward,
from the outer boundary circle to the inner one, and mark the non-empty
union of curves $\pi^{-1}(C_x)\cap q(Q)$.

This picture is called the \textit{portrait of $x$}. Two portraits are
considered equal if they differ by an orientation-preserving fiberwise
diffeomorphism of the torus.

\begin{figure}
		\centering
		\includegraphics[width=0.96\textwidth]{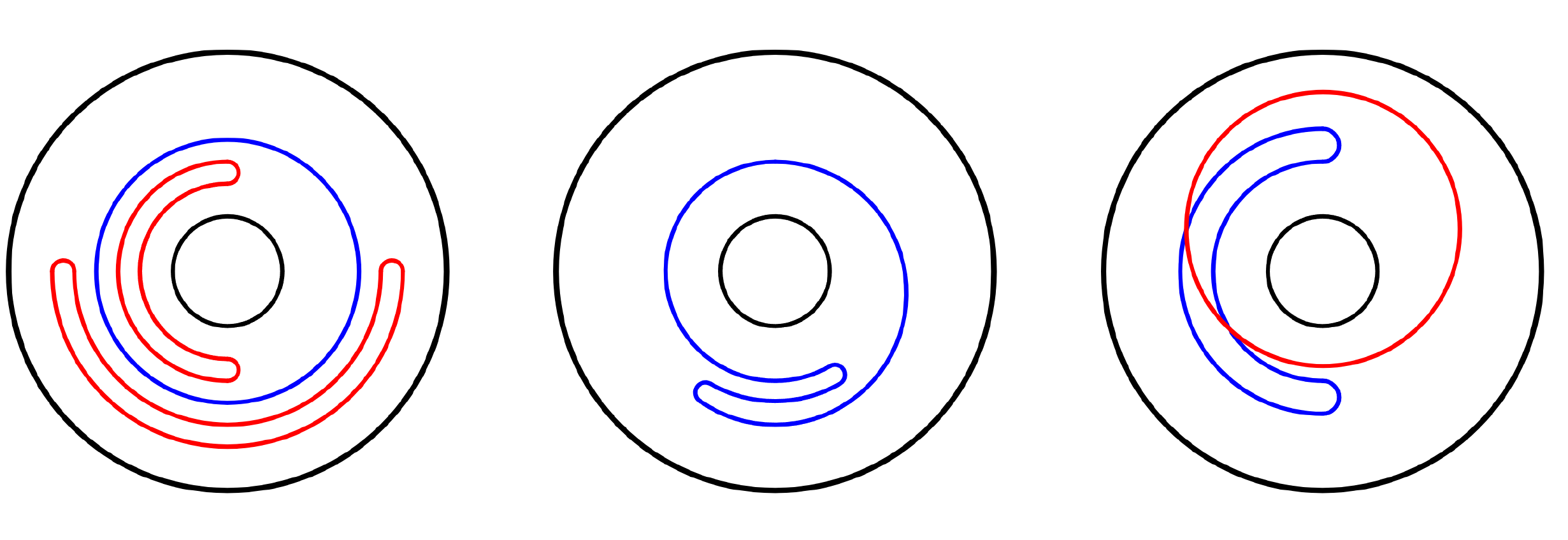}
		\caption{The non-symmetric portraits. Left: type $I$, with two
fold curves and one simple circle; at least one simple circle must be
present. Center: type $II$, with a pleat and no simple circles. Right:
type $III$, with a simple circle crossing a fold curve.}
		\label{fig:nonsymm}
	\end{figure}

	\begin{lemma}\label{LemmaPortr}
		\begin{enumerate}
			\item The portrait of an intersection of two projected folds is
as in the left-hand part of Figure~\ref{fig:nonsymm}, up to adding or
removing simple circles. At least one simple circle is present because
$\pi\circ q(Q)=B$.

Such a singular vertex is called \textbf{type $I$}.
			
			\item The portrait of a projected pleat is, up to reflection,
as in the center of Figure~\ref{fig:nonsymm}. The number of simple
circles may vary.

Such a vertex is called \textbf{type $II$, right} if its portrait is
the one shown in the figure, and \textbf{type $II$, left} if its
portrait is its reflection.

			\item The portrait of the projection of an intersection between
a fold and a regular sheet is, up to reflection, as in the right-hand
part of Figure~\ref{fig:nonsymm}.
			
			Such a vertex is called \textbf{type $III$, right} if its
portrait is the one shown in the figure, and \textbf{type $III$, left}
if its portrait is its reflection, in both cases up to adding simple
circles.

		\end{enumerate}
		
	\end{lemma}

	Singular vertices of types $I$, $II$, and $III$ are called
\textit{essential vertices}.

\subsection{Closed quasisection and the local formula}
\label{SecWeightsClosed}

In this subsection $\mathcal Q=(Q,q)$ is again a generic closed
quasisection.
	
	Let $x$ be a singular vertex of type $I$.

	Move a point $y$ counterclockwise along $C_x$, starting in a sector
with no projected fold curves, equivalently, in a sector where
$\pi^{-1}(y)\cap q(Q)$ has the minimum number of points. Order the two
fold curves according to when they are first encountered by $y$.
	
Let $n(x)$ be the number of simple circles between the first and second
fold curves, counted from the first fold curve in the direction of the
fiber. Similarly, let $k(x)$ be the number of simple circles between
the two fold curves, counted from the second fold curve in the direction
of the fiber.
	
\medskip

	\begin{lemma}\begin{enumerate}
		               \item The portrait of an essential vertex of type $I$ is
determined by the two non-negative integers $n$ and $k$ defined above. Their sum $n+k$ is non-zero.
		
		               \item The portrait of an essential vertex of type
$II$ is determined by, and determines, its side (right or left) and the
number $r\geq0$ of simple circles.
		               \item The portrait of an essential vertex of type
$III$ is determined by, and determines, its side (right or left) and the
number $r\geq0$ of simple circles.
	             \end{enumerate}
\end{lemma}

\medskip

In \cite{PanShamTur}, weights are assigned to the essential singular vertices:

	\begin{definition}\begin{enumerate}
			\item For an essential vertex $x$ of type $I$, set $$\mathcal{W}_{ff}(x)=\mathcal{W}_{ff}(n,k)=\frac{4(n-k)}{(n+k)(n+k+2)(n+k+4)}.$$
			
			\item For a right essential vertex $x$ of type $II$, set $$\mathcal{W}_{p}(x)=\mathcal{W}_{p,R}(r)=\frac{2}{(r+1)(r+3)}.$$
			
			For a left essential vertex $x$ of type $II$, set $$\mathcal{W}_{p}(x)=\mathcal{W}_{p,L}(r)=\frac{-2}{(r+1)(r+3)}.$$
			
			\item For a right essential vertex $x$ of type $III$, set $$\mathcal{W}_{fs}(x)=\mathcal{W}_{fs}(r,R)= \frac{2}{(r+1)(r+3)}.$$
			
			For a left essential vertex $x$ of type $III$, set $$\mathcal{W}_{fs}(x)=\mathcal{W}_{fs}(r,L)=\frac{-2}{(r+1)(r+3)}.$$
			
		\end{enumerate}
		
	\end{definition}

	\begin{theorem}[The local formula for closed quasisections, \cite{PanShamTur}]\label{ThmMain}
		The Euler number of a circle bundle with a generic closed
quasisection $\mathcal{Q}$ is equal to the sum of the weights of its
essential singular vertices:
		
		$$\mathcal{E}= \sum_I \mathcal{W}_{ff}(x_i)+ \sum_{II} \mathcal{W}_{p}(x_i)+\sum_{III} \mathcal{W}_{fs}(x_i).$$
		
		\end{theorem}

\subsection{Bordered quasisections and the local formula}\label{SecBordered}

Now $\mathcal Q=(Q,q)$ is a generic quasisection whose boundary may be non-empty. In this setting the
singular curve contains also the projection of $q(\partial Q)$. The definitions of the three
interior types $I$, $II$, and $III$ remain unchanged.
Four additional types of essential vertices can occur:
\begin{enumerate}
    \item $I_{bb}$: a transverse intersection of two projected boundary
    curves;
    \item $I_{bf}$: a transverse intersection of a projected boundary
    curve and a projected fold;
    \item $III_{bs}$: the projection of an intersection of a boundary
    curve with a regular sheet;
    \item $IV_{tf}$: the projection of a point at which a fold
    terminates on the boundary of $Q$ (a $B_2$-point).
\end{enumerate}
The types $I_{bf}$, $III_{bs}$, and $IV_{tf}$ have right and left
versions, distinguished by reflection of their local portraits. For
$I_{bb}$ and $I_{bf}$, the numbers $n$ and $k$ count the simple circles
between the two ordered curves in the two directions along the fiber,
as for type $I$; in particular, $n+k>0$. For types $III_{bs}$ and $IV_{tf}$,
$r\geq 0$ is the number of simple circles in the portrait. The
right--left conventions are those of \cite[Section 6]{PanShamTur}.

The weights of the new types of singular vertices are
\begin{align*}
\mathcal{W}_{bb}(n,k)
    &=\frac{n-k}{(n+k)(n+k+1)(n+k+2)},\\
\mathcal{W}_{bf}(n,k,R)
    &=\frac{n-k}{(n+k)(n+k+3)}
      \left(\frac{1}{n+k+1}+\frac{1}{n+k+2}\right),\\
\mathcal{W}_{bf}(n,k,L)&=-\mathcal{W}_{bf}(n,k,R),\\
\mathcal{W}_{bs}(r,R)&=\frac{1}{(r+1)(r+2)},\\
\mathcal{W}_{bs}(r,L)&=-\mathcal{W}_{bs}(r,R),\\
\mathcal{W}_{tf}(r,R)&=\frac{1}{2(r+1)(r+2)},\\
\mathcal{W}_{tf}(r,L)&=-\mathcal{W}_{tf}(r,R).
\end{align*}
The weights of types $I$, $II$, and $III$ remain unchanged.

\begin{theorem}[The local formula for possibly bordered quasisections, \cite{PanShamTur}]\label{ThmBorderedLCF}
\label{ThmBorderedFormula}
For a generic bordered quasisection $\mathcal{Q}$ of $E\rightarrow B$,
the Euler number is equal to the sum of weights:
\begin{align*}
\mathcal{E}={}&
  \sum_{I}\mathcal{W}_{ff}(x_i)
 +\sum_{II}\mathcal{W}_{p}(x_i)
 +\sum_{III}\mathcal{W}_{fs}(x_i)\\
 &+\sum_{I_{bb}}\mathcal{W}_{bb}(x_i)
 +\sum_{I_{bf}}\mathcal{W}_{bf}(x_i)
 +\sum_{III_{bs}}\mathcal{W}_{bs}(x_i)
 +\sum_{IV_{tf}}\mathcal{W}_{tf}(x_i).
\end{align*}
Clearly, if $\partial Q=\varnothing$,
the four boundary terms vanish and the formula reduces to
Theorem~\ref{ThmMain}.
\end{theorem}

\section{Proof of Theorem~\ref{ThmMinimal}: Lower bounds for the number of singular vertices}\label{SecLowerBounds}

Let $\mathcal{Q}$ be a generic (possibly bordered) quasisection.
The definition of the weights implies immediately:
\[
\begin{aligned}
 |\mathcal{W}_{ff}|&\leq\frac4{15}, &
 |\mathcal{W}_{p}|&\leq\frac23, &
 |\mathcal{W}_{fs}|&\leq\frac23, &
 |\mathcal{W}_{bb}|&\leq\frac16,\\
 |\mathcal{W}_{bf}|&\leq\frac5{24}, &
 |\mathcal{W}_{bs}|&\leq\frac12, &
 |\mathcal{W}_{tf}|&\leq\frac14.
\end{aligned}
\]
Consequently, the weight of every singular vertex $x$ satisfies
\[
 |\mathcal{W}(x)|\leq\frac23,
\]
and we have
\[
 |\mathcal{E}|
 \leq\sum_x|\mathcal{W}(x)|
 \leq\frac23N(\mathcal{Q}).
\]
Since $N(\mathcal{Q})$ is an integer, we have
\begin{equation}\label{EqBorderedLowerBound}
 \mathcal{N}_{\mathrm{b}}(E\rightarrow B)\geq
 \left\lceil\frac{3|\mathcal{E}|}{2}\right\rceil,
\end{equation}

The same is true for closed quasisections:
\begin{equation}\label{EqClosedEvenLowerBound}
 \mathcal{N}_{\mathrm{cl}}(E\rightarrow B)
 \geq\left\lceil\frac{3|\mathcal{E}|}{2}\right\rceil,
\end{equation}

but this is not yet the desired formula from item~\ref{ThmMainClosed} of Theorem~\ref{ThmMinimal} for odd $\mathcal{E}$, so we need a further analysis.

\medskip

Assuming that $Q$ is closed, recall that the degree of the
quasisection, denoted by $\deg(\mathcal{Q})$, is either the integral
degree or the degree modulo~$2$ of the map
$\pi\circ q\colon Q\to B$.

It is proven in
\cite[Theorem~12]{PanShamTur} that
\begin{enumerate}
     \item if $Q$ is orientable and
$\mathcal E\ne0$, the integral degree vanishes;
     \item if $Q$ is
non-orientable and $\mathcal E$ is odd, its degree modulo~$2$ is zero.
   \end{enumerate}

If $\mathcal{E}$ is odd, then for
each regular point $x$ of the base, the number of points in the preimage $\bigl|(\pi\circ q)^{-1}(x)\bigr|$ is even, and therefore all essential vertices have absolute weight at most
$1/4$. This follows from the definitions of the weights and from
\[
 r\geq1\quad\text{for types $II$ and $III$},
 \qquad n+k\geq2\quad\text{for type $I$}.
\]

Thus the local formula gives $N(\mathcal Q)\geq4|\mathcal E|$ for
every closed quasisection when $\mathcal E$ is odd, that is,
\begin{equation}\label{EqClosedOddLowerBound}
 \mathcal{N}_{\mathrm{cl}}(E\rightarrow B)
 \geq4|\mathcal{E}|.
\end{equation}

Thus we have proved all the lower bounds in
Theorem~\ref{ThmMinimal}.

\section{Proof of Theorem~\ref{ThmMinimal}: Examples realizing the equalities}
It remains to construct examples of quasisections with the desired numbers of singular vertices.
We first construct examples for Euler numbers $2$ and $1$, and then
combine these basic examples for a bundle with arbitrary
Euler number.

\subsection{A quasisection with the smallest possible number of singular vertices for $\mathcal{E}=2$}\label{SecExampleTwo}

Assume that $\mathcal{E}=2$. Choose a small\footnote{Although no metric on $B$ is fixed, one may imagine choosing an auxiliary metric and taking $D$ to be a small disk with respect to it.} disk in the base
$D\subset B$. Take two trivializations of the bundle, one over $D$ and the other one over
$B\setminus\operatorname{Int}D$.

Each of the trivializations gives a diffeomorphism of the gluing torus
\[
 T=\pi^{-1}(\partial D)=\partial D\times S^1
\]
with the standard torus $S^1\times S^1$. The restrictions of the two trivializations to the boundary $\partial D$ differ (up to homotopy) by the square of the Dehn
twist with respect to the fiber.

Now choose a constant section $s$ over $B\setminus\operatorname{Int}D$. Its
 graph is some bordered surface $S\subset E$ whose
boundary $\partial S$ lies on the gluing torus $T$.

\begin{figure}[!b]
\centering
\includegraphics[width=0.48\textwidth]{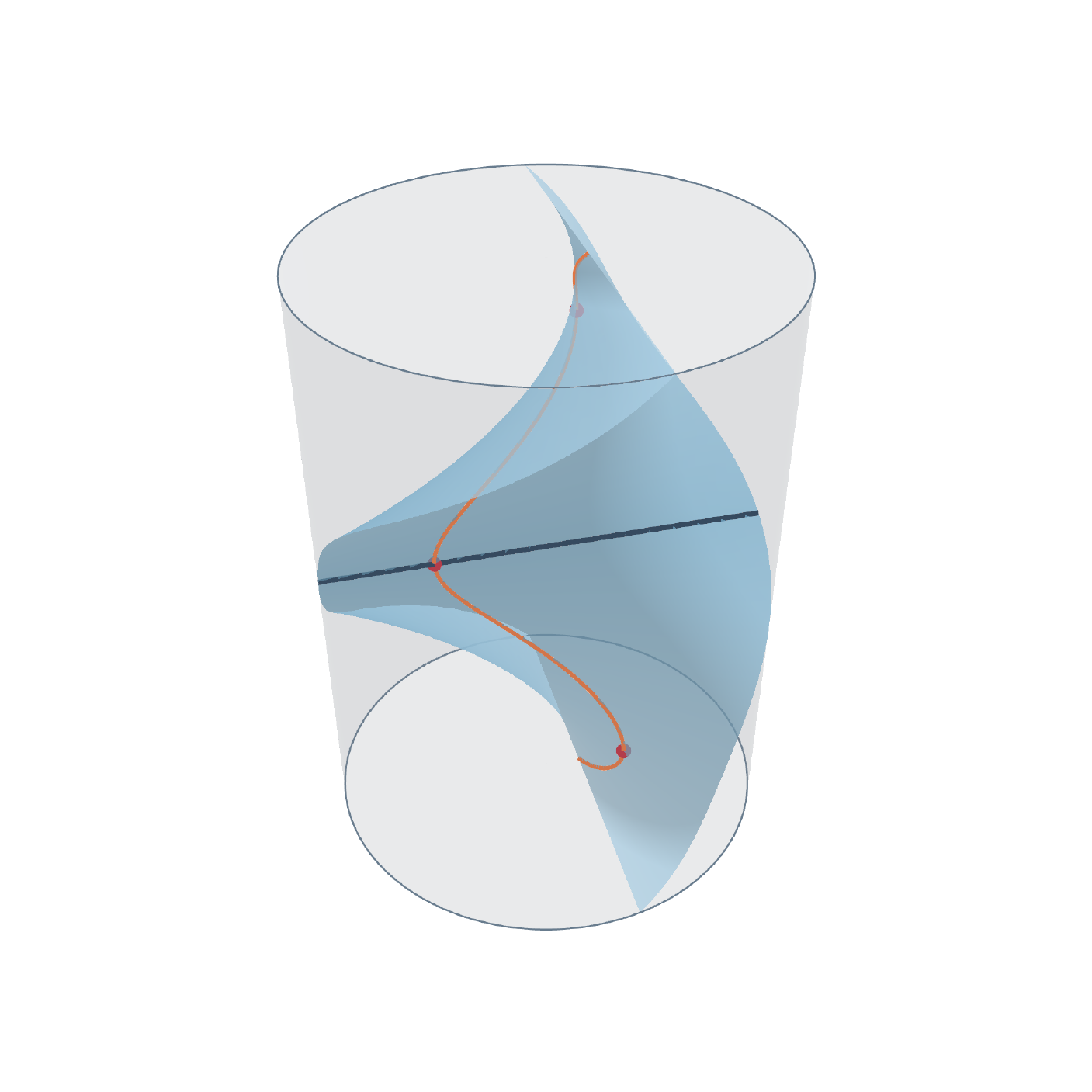}\hfill
\includegraphics[width=0.48\textwidth]{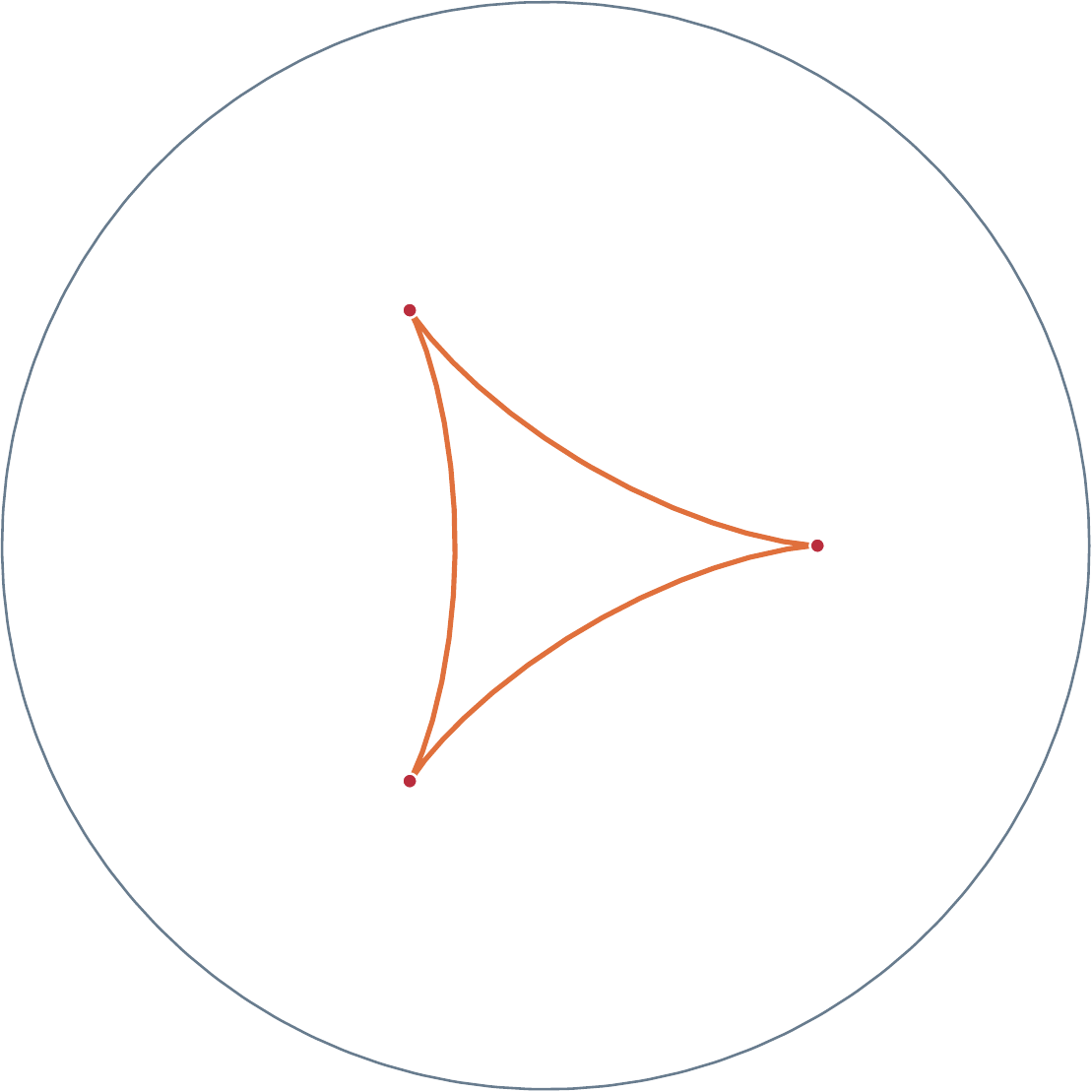}
\caption{The  construction for $\mathcal{E}=2$. Left: the
smooth M\"obius band, its fold curve, and the ambient solid torus $\pi^{-1}(D)$. Right:
the projection of the fold curve to $D$; its three cusps are the three
vertices of type $II$.}
\label{FigMobiusBlock}
\end{figure}

Let
$a=[\partial D\times\{*\}]$ and $b=[\{*\}\times S^1]$ be the standard generators of $T$. In the
coordinates induced by the trivialization over $D$, the curve
$\partial S$ is homotopic to $a+2b$. In Figure~\ref{FigMobiusBlock}, it
is visualized as a curve spiraling on $T$. To build a quasisection, we construct inside
$D\times S^1$ a surface $S'$ whose boundary equals $\partial S$. Then gluing $S$ and $S'$ along their boundaries yields a closed quasisection.

Identify $D$ with the unit disk and write the circle
fibers as $\mathbb R/2\pi\mathbb Z$. So a point of
$\pi^{-1}(D)=D\times S^1$ has coordinates $(x,\alpha)$, where
$|x|\leq1$ and $\alpha$ is taken modulo $2\pi$.

Consider the M\"obius band
\[
 M=\bigl([0,1]\times[-1,1]\bigr)/\bigl((1,-u)\sim(0,u)\bigr).
\]
For $t\in[0,1]$, set
\[
 \begin{aligned}
  \theta&=\frac{\pi}{3}+\pi t,
  &\mathbf n&=(\cos\theta,\sin\theta),\\
  \mathbf n^\bot&=(-\sin\theta,\cos\theta),
  &h&=-\frac16\cos3\theta.
 \end{aligned}
\]
Set also
\begin{equation}\label{EqMobiusBlock}
 \mathcal S:M\longrightarrow D\times S^1,\qquad
 \mathcal S(t,u)=
 \left(h\mathbf n+u\sqrt{1-h^2}\,\mathbf n^\bot,
 e^{2\pi it}\right).
\end{equation}

Let us comment on this formula. 
For fixed $t$, the map $h\mathbf n+u\sqrt{1-h^2}$
parameterizes a chord of the disk $D$. The unit vector $\mathbf n$ is the
normal to this chord, \(h\) is the signed distance of the chord from the origin, and
$\mathbf n^\bot$ gives its direction. The value $u=0$ corresponds to its midpoint, while
$u=-1$ and $u=1$ correspond to its endpoints. The endpoints lie on $\partial D$.

The image $S'=\operatorname{Im}\mathcal S$ is a smooth M\"obius band in the
solid torus $D\times S^1$. Geometrically, imagine a moving horizontal segment in
the solid torus whose projection is a chord of $D$ tangent to the
three-cusped curve shown in Figure~\ref{FigMobiusBlock}. The union of
these segments is $S'$. The homotopy class of its boundary is $a+2b$,
hence $S'$ can be glued to $S$ along their boundaries.

We arrive at a closed quasisection.
To find the critical curve, write the first component of $\mathcal{S}$ as
\[
 X(\theta,u)=h\mathbf{n}+u\sqrt{1-h^2}\,\mathbf{n}^\bot.
\]
Since $(\mathbf{n})'=\mathbf{n}^\bot$ and $(\mathbf{n}^\bot)'=-\mathbf{n}$, we have
\[
 X_u=\sqrt{1-h^2}\,\mathbf{n}^\bot,
 \qquad
 X_\theta=\bigl(h'-u\sqrt{1-h^2}\bigr)\mathbf{n}
 +\left(h+u\frac{d}{d\theta}\sqrt{1-h^2}\right)\mathbf{n}^\bot.
\]
Consequently, $X_u$ and $X_\theta$ are linearly dependent exactly when
$u=h'/\sqrt{1-h^2}$. Since $|h'|\leq1/2$ and
$\sqrt{1-h^2}\geq\sqrt{35}/6$, this value lies
in the interior of the chord. Substitution gives the critical curve
\[
 \gamma(\theta)=X\left(\theta,\frac{h'}{\sqrt{1-h^2}}\right)
 =h\mathbf{n}+h'\mathbf{n}^\bot.
\]
Since
\[
 \gamma'(\theta)=(h+h'')\mathbf{n}^\bot=\frac43\cos3\theta\,\mathbf{n}^\bot,
\]
and $\theta$ ranges from $\pi/3$ to $4\pi/3$, the critical curve has exactly
three cusps that correspond to three pleats. They are singular vertices of type $II$
with $r=0$, have the same sign, and each has weight $2/3$. The glued
surface $S\bigcup S'$ is therefore a closed quasisection of $E\rightarrow B$ with
exactly three singular vertices; see Figure~\ref{FigMobiusBlock}.

Clearly, the same number of vertices is realized in the bordered class. For instance, one
may add an embedded regular disk whose image lies over a region already
covered by the closed quasisection. This creates a non-empty boundary
but no new singular vertices.

\subsection{A quasisection with the smallest possible number of singular vertices for $\mathcal{E}=1$}\label{SecExampleOne}

Assume that $\mathcal{E}=1$. Again, choose a
small disk $D\subset B$, set $B_0=B\setminus\operatorname{Int}D$, and
take trivializations of the bundle over $D$ and over $B_0$. This time
their restrictions to $\partial D$ differ (up to homotopy) by one
positive Dehn twist with respect to the fiber.

\begin{figure}[!t]
\centering
\includegraphics[width=0.48\textwidth]{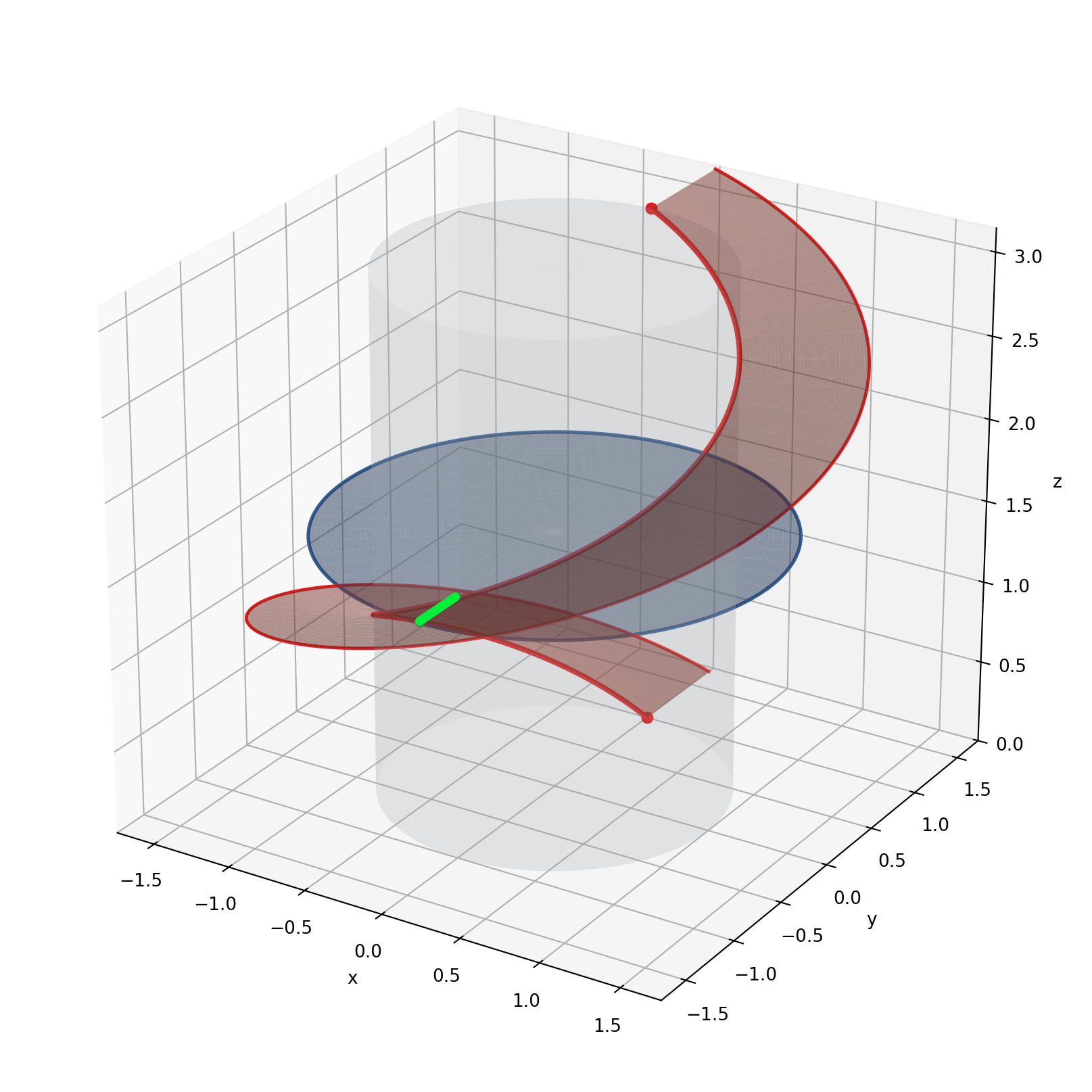}\hfill
\includegraphics[width=0.48\textwidth]{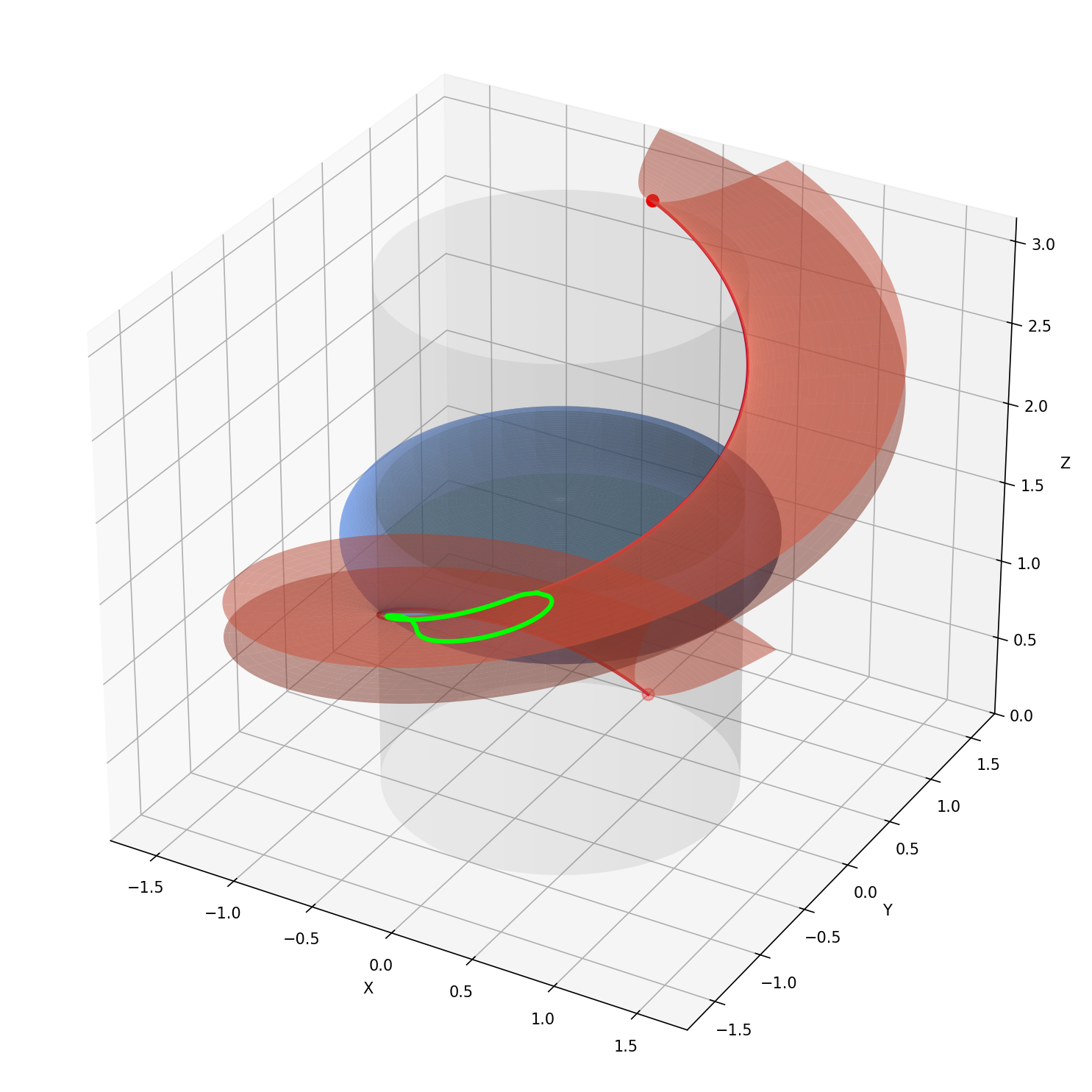}
\caption{Local constructions for $\mathcal{E}=1$. Left: the bordered
construction formed by two disks and the normal strip produced by the
twist. Right: the closed construction formed by the boundaries of two
tubular neighborhoods. The gluing torus and the relevant intersections
are indicated in both pictures, which are shown in the coordinates
induced by the trivialization over $D$.}
\label{FigEulerOneBlocks}
\end{figure}

We use the same notation as above:
\[
 T=\pi^{-1}(\partial D),\qquad
 a=[\partial D\times\{*\}],\qquad b=[\{*\}\times S^1]
\]. 
Choose a constant section over $B_0$ and denote its graph by
$S$. In the coordinates induced by the trivialization over $D$, its
boundary is homotopic to $a+b$, since the Dehn twist adds one copy of
the fiber class $b$ to $a$.

Let $V=\pi^{-1}(D)$ and choose a constant section disk
$S'\subset V$ over $D$. Its boundary represents the class $a$.
Consequently, $\partial S'$ and $\partial S$ have algebraic
intersection one. Without loss of generality, we choose these curves
to intersect transversely at exactly one point.

Extend $S'$ by attaching a narrow collar across $T$, and modify the surface near the attaching curve so that the resulting disk, still denoted by $S'$, is smooth and
satisfies
\[
 \pi(\partial S')\subset\operatorname{Int}B_0.
\]
The collar may be chosen in general position so that $S\cap S'$ is one
transverse arc joining $\partial S$ to $\partial S'$. Its two
endpoints project to two right vertices of type $III_{bs}$ with $r=0$.
This construction is shown on the left-hand side of
Figure~\ref{FigEulerOneBlocks}.

The construction has only two singular vertices, each of weight
\[
 \mathcal{W}_{bs}(0,R)=\frac12,
\]
so their total weight is $1$. This gives the required example of a bordered quasisection with the smallest possible number of singular vertices for $\mathcal{E}=1$.

For the closed construction, start with the two embedded surfaces
$S$ and $S'$ from the bordered construction and take the
boundaries of sufficiently small tubular neighborhoods of the two
surfaces separately. After a small generic
perturbation, the two resulting closed surfaces form a closed
quasisection with exactly four singular vertices, all of type $III$
with $r=1$ and weight $1/4$. This construction is shown on the
right-hand side of Figure~\ref{FigEulerOneBlocks}.

\begin{figure}[H]
\centering
\includegraphics[width=0.96\textwidth]{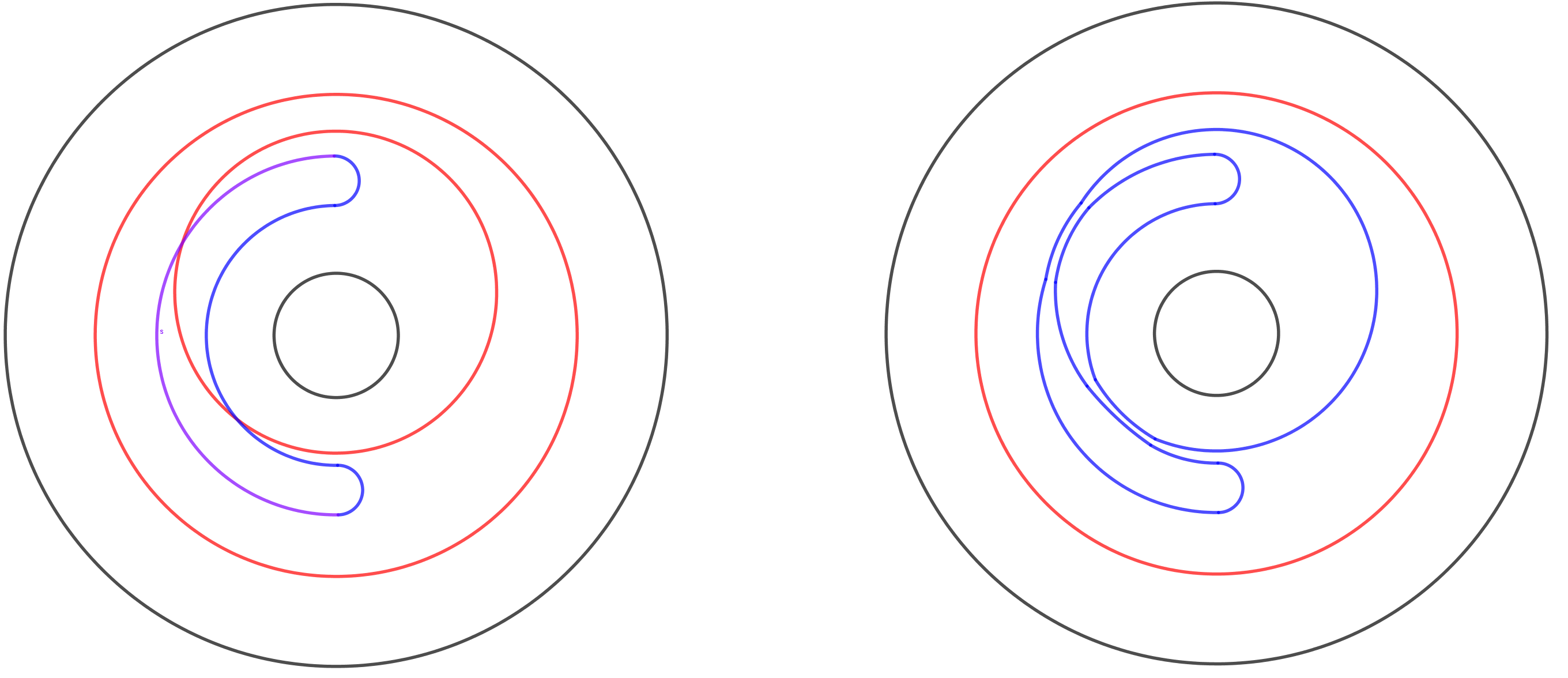}
\caption{Eliminating self-intersections by surgery, as seen in the
portraits of singularities. Left: before the surgery. Right: after the
surgery.}
\label{FigSurgery}
\end{figure}

One can go even further and, using an appropriate surgery, eliminate
the self-intersections of the quasisection. The surgery can be
visualized through the portraits of singularities; see
Figure~\ref{FigSurgery}.

\subsection{The general case}\label{SecExamplesArbitrary}

Choose pairwise disjoint disks in $B$. The bundle is trivial over their
complement, and its Euler number is the sum of the degrees of the
clutching twists around the disks. We may therefore insert independent
copies of the $\mathcal{E}=1$- and $\mathcal{E}=2$- constructions from
Sections~\ref{SecExampleTwo} and~\ref{SecExampleOne} and attach them to a common partial section over
the complement. Since their supports are disjoint, the numbers of their
singular vertices add up. For negative $\mathcal E$, use negative powers or simply reverse the orientation of the base.

For bordered quasisections, write $|\mathcal E|=2m$ or $2m+1$. Use $m$
copies of the construction in Section~\ref{SecExampleTwo} and, in the odd case, one
bordered copy of the construction in Section~\ref{SecExampleOne}. This gives $3m$ or
\[
 3m+2=\left\lceil\frac{3|\mathcal{E}|}{2}\right\rceil
\]
vertices, respectively. For $\mathcal{E}=0$, a section together with an embedded
regular disk is a bordered quasisection without singular vertices.
Together with~\eqref{EqBorderedLowerBound}, these constructions prove
item~\ref{ThmMainBordered} of Theorem~\ref{ThmMinimal}.

For closed quasisections with $|\mathcal E|=2m$, use $m$ copies of the
closed construction in Section~\ref{SecExampleTwo}; the result has $3m$ vertices. If
$|\mathcal E|$ is odd, use $|\mathcal E|$ disjoint copies of the closed
Euler-number-one construction from Section~\ref{SecExampleOne}; the result has
$4|\mathcal E|$ vertices. Finally, for $\mathcal E=0$ a global section has
no singular vertices. These examples attain
\eqref{EqClosedEvenLowerBound} and~\eqref{EqClosedOddLowerBound}, and
complete the proof of Theorem~\ref{ThmMinimal}.

\textbf{Acknowledgement.}
This work is supported by the Russian Science Foundation (project 25-11-00058).

\medskip

\textbf{AI Usage Declaration .}
We hereby declare that no AI was used in the completion of this research.

\end{document}